\documentclass[leqno,11pt]{amsart}
\usepackage[T1]{fontenc}
\usepackage{lmodern}
\usepackage[margin=25mm]{geometry}
\usepackage{amsmath,amssymb,amsthm,mathtools}
\usepackage{microtype}
\usepackage{xcolor}
\usepackage{enumitem}
\usepackage{fancyhdr}
\usepackage[hidelinks]{hyperref}
\hypersetup{pdftitle={A double-logarithmic lifespan bound for axisymmetric Euler equations with a small initial swirl}}
\numberwithin{equation}{section}
\newtheorem{theorem}{Theorem}[section]
\newtheorem{lemma}[theorem]{Lemma}

\newtheorem{corollary}[theorem]{Corollary}
\theoremstyle{remark}
\newtheorem{remark}[theorem]{Remark}

\newcommand{\norm}[2]{\left\|#1\right\|_{#2}}

\newcommand{\er}{\bl{e}_r}
\newcommand{\etheta}{\bl{e}_\th}
\newcommand{\ez}{\bl{e}_z}

\newcommand{\R}{\mathbb R}

\newcommand{\dy}{\,{\rm d}y}
\newcommand{\dx}{\,{\rm d}x}

\def\bR{{\mathbb R}}

\def\al{\alpha}

\def\ve{\varepsilon}

\def\th{\theta}

\def\lt{\left}
\def\rt{\right}

\def\i{\infty}
\def\p{\partial}
\def\f{\frac}
\def\na{\nabla}

\def\O{\Omega}
\def\s{\sqrt}
\def\q{\quad}
\def\qq{\qquad}

\def\bl{\boldsymbol}

\def\mK{\mathcal{K}}

\def\r{\rho}

\allowdisplaybreaks
\def\cd{\cdot}

\def\bm{\boldsymbol}

\def\dx{\,\mathrm{d}x}
\def\d{\,\mathrm{d}}
\def\ba{\begin{equation}\begin{aligned}}
	\def\ea{\end{aligned}\end{equation}}
\def\bn{$$\begin{aligned}}
\def\en{\end{aligned}$$}
\def\l{\label}
\def\ed{\buildrel\hbox{\footnotesize def}\over =}
\def\ef{\eqref}

\def\be{\begin{equation}}
\def\ee{\end{equation}}
\def\wt{\widetilde}
\newcommand{\beq}{\begin{equation}}
	\newcommand{\eeq}{\end{equation}}
\newcommand{\ben}{\begin{eqnarray}}
	\newcommand{\een}{\end{eqnarray}}
\newcommand{\beno}{\begin{eqnarray*}}
	\newcommand{\eeno}{\end{eqnarray*}}

\renewcommand{\div}{\mathop{\rm div}\nolimits}

\renewcommand{\geq}{\geqslant}
\renewcommand{\d}{\,{\rm d}}

\def\XXint#1#2#3{{\setbox0=\hbox{$#1{#2#3}{\int}$}	\vcenter{\hbox{$#2#3$}}\kern-.5\wd0}}

\title{A double-logarithmic lifespan bound for axisymmetric Euler equations with a small initial swirl}
\author[Z. Li]{Zijin Li}
\address[Z. Li]
{School of Mathematics and Statistics, Nanjing University of Information Science and Technology, Nanjing
	210044, China}
\email{zijinli@nuist.edu.cn}

\author[X. Xiao]{Xu Xiao}
\address[X. Xiao]
{School of Mathematics and Statistics,
Nanjing University of Information Science and Technology,
Nanjing 210044, China}
\email{xiaoxu@nuist.edu.cn}

\author[L. Yang]{Linbin Yang}
\address[L. Yang]
{School of Mathematics and Statistics,
Nanjing University of Information Science and Technology,
Nanjing 210044, China}
\email{yanglinbin@nuist.edu.cn}

\date{}
\begin{document}
\maketitle
\vspace{-1.4em}
\begin{abstract}
We study the axisymmetric incompressible Euler equations in $\mathbb{R}^3$ with a small initial swirl. Fix $s>7/2$, and let $\boldsymbol{v}_0\in H^s(\mathbb{R}^3)$ be axisymmetric and divergence-free. If $\varepsilon=\|{\nabla\times(v_{0,\theta}\boldsymbol{e}_{\theta})}\|_{L^\infty}$ is sufficiently small, then the lifespan $T_*$ of the corresponding $H^s$ solution satisfies $$ T_{*}\geq\frac{C_s}{1+\|{\boldsymbol{v}_0}\|_{H^s}} \left[ \log\log (\varepsilon^{-1}) \right]^{3/7}\,. $$ This removes the assumption by Danchin (Proc. Amer. Math. Soc. 141:1979–1993, 2013) that the domain must be distant from the axis of symmetry. Moreover, even for $\boldsymbol{v}_0\in H^3$, our result improves the fourfold iterated-logarithmic lower bound for the lifespan obtained by Li and Zhou (Z. Angew. Math. Phys. 75:219, 2024) to a threefold iterated-logarithmic bound.
\end{abstract}

\tableofcontents

\section{Introduction and main result}

In this paper, we consider the incompressible Euler equations
\be
\l{eq:euler}
\begin{cases}
\p_t \bl{v}+\bl{v}\cd\na \bl{v}+\na p=0,\\
\div \bl{v}=0,\\
\bl{v}(0,x)=\bl{v}_0,
\end{cases}
\qq (t,x)\in[0,T)\times\R^3,
\ee
where $\bl{v}$ is the velocity and $p$ is the pressure. For
divergence-free $\bl{v}_0\in H^s(\R^3)$, $s>5/2$, the classical
local theory gives a unique maximal solution
$$
\bl{v}\in C([0,T_{*});H^s(\R^3))
\cap C^1([0,T_{*});H^{s-1}(\R^3))\,.
$$
see \cite{KP} for the local theory in Sobolev spaces of real order and \cite[Section~3.2]{MB} for the energy method.

We call a solution to \eqref{eq:euler} is axially symmetric, if and only if the solution $\bl{v}$ enjoys the following expressing:
$$
\bl{v}=v_r(t,r,z)\er+v_\th(t,r,z)\etheta+v_z(t,r,z)\ez.
$$
Here $v_r,v_\th, v_z$ are scalar functions of $(t,r,z)$, where $(r,\th,z)$ is the cylindrical coordinates:
$$
r=\sqrt{x_1^2+x_2^2}, \qq \th=\arctan\f{x_2}{x_1}, \qq z=x_3,
$$
and the basis vectors $\er,\etheta,\ez$ are
$$
\er=\lt(\f{x_1}{r},\f{x_2}{r},0\rt),\qq
\etheta=\lt(-\f{x_2}{r},\f{x_1}{r},0\rt),\qq
\ez=(0,0,1).
$$
In this cylindrical coordinate system, system \eqref{eq:euler} can be rewriting into
\be
\l{eq:cylindrical}
\left\{\begin{aligned}
&\p_t v_r+(v_r\p_r+v_z\p_z)v_r
-\f{v_\th^2}{r}+\p_r p=0,\\
&\p_t v_\th+(v_r\p_r+v_z\p_z)v_\th
+\f{v_rv_\th}{r}=0,\\
&\p_t v_z+(v_r\p_r+v_z\p_z)v_z+\p_z p=0,\\
&\p_r v_r+\f{v_r}{r}+\p_z v_z=0.
\end{aligned}\right.
\ee

Below, we denote
$$
\bl{b}=v_r\er+v_z\ez,\qq
\bl{w}=\na\times \bl{v}=w_r\er+w_\th\etheta+w_z\ez.
$$
Here
$$
w_r=-\p_z v_\th,\qq w_\th=\p_z v_r-\p_r v_z,\qq w_z=\p_r v_\th+\f{v_\th}{r}.
$$
We also write
$$
\O=\f{w_\th}{r},\qq J=\f{w_r}{r},\qq r>0,
$$
and denote their initial values by
$$
\O_0=\f{w_{0,\th}}{r},\qq J_0=\f{w_{0,r}}{r}.
$$
By incompressibility, the field $\bl{b}$ is divergence free. For a flow without swirl, the vorticity equation gives
$$
(\p_t+\bl{b}\cd\na)\O=0.
$$
This transport equation is central to the global regularity theory for axisymmetric Euler flows without swirl. The classical results of Ukhovskii and Yudovich \cite{UY} and Shirota and Yanagisawa \cite{SY} give global solutions under suitable regularity assumptions; see also \cite{MB}. Abidi, Hmidi and Keraani \cite{AHK} extended the theory to critical Besov spaces under an additional Lorentz-space assumption on $\O_0$. For smooth solutions with swirl, global regularity in $\R^3$remains open. A natural question is therefore how long a solution exists when its initial swirl is small, while its meridional velocity need not be small.

Danchin \cite{Danchin} studied this question on bounded axisymmetric domains separated from the symmetry axis. Under suitable Besov smallness assumptions on the swirl, he obtained a double-logarithmic lower bound for the lifespan. Related long-time existence questions also arise for perturbations of the two-dimensional Euler equations. Later Subsequently, Li and Zhou \cite{LZ} extended this result to the whole space, thereby removing the restriction that the support of the velocity field must remain away from the axis of symmetry: Suppose an axisymmetric divergence-free velocity field $\bm{v}_0\in H^3(\R^3)$ with
$$
\norm{\na\times(v_{0,\th}\etheta)}{L^\i}=\ve\ll1,
$$
then the lifespan lower bound satisfies
$$
T_{*}\ge\f{C}{1+\norm{\bl{v}_0}{H^3}}
\log\log\log\log(\ve^{-1})\,.
$$
Recently, Bae, Lopes Filho, Mazzucato and Nussenzveig Lopes \cite{BLMN} gave a proof of long-time existence for the ideal Boussinesq and density-dependent Euler equations with small temperature or density perturbations, without requiring the initial velocity to be small.

The lifespan problem considered in this paper is closely related to the long-time growth of axisymmetric solutions to the Euler equations without swirl. Although axisymmetric Euler flows without swirl are globally regular under the assumptions stated above, $\norm{w_\th(t,\cd)}{L^\i}$ need not remain uniformly bounded in time. Indeed, transport of $\O$ allows vorticity amplification when fluid particles move away from the symmetry axis. Childress \cite{Childress} conjectured that the maximal vorticity growth rate is $t^{4/3}$. The eroding vortex dipole model and numerical simulations in \cite{CGV} later provided evidence that this rate can be attained.

Recently, Lim and Jeong \cite{LJ} proved
$$
\norm{w_\th(t,\cd)}{L^\i}\le C(1+t)^{4/3}
$$
for axisymmetric flows without swirl whose initial vorticity is compactly supported and satisfies
$\O_0\in L^\i(\R^3)$.
The key step in their proof is to use the conserved $L^2$ norm of $\bl{v}$ to estimate $v_r$, which in turn controls the growth of $\norm{w_\th(t,\cd)}{L^\i}$. Shao, Wei and Zhang \cite{SWZ} obtained the same upper bound without the compact-support assumption, for $\bl{v}_0\in L^2(\R^3)\cap C^{1,\al}(\R^3)$ with $0<\al<1$, assuming that both $w_{0,\th}$ and $\O_0$ belong to $L^\i(\R^3)$. It remains open whether the $t^{4/3}$ upper bound is optimal.

Inspired by these results, for $\bl{v}_0\in H^s(\R^3), s>7/2$, we improve the lifespan estimate for small initial swirl by establishing an improved estimate of $\int_0^t\norm{w_\th(t^\prime,\cd)}{L^\i}\d t^\prime$. The main difficulty is that the equation for $\O$ contains a source term involving the swirl.
To control this source term, we estimate $\norm{\O}{L^\i}$, $\norm{J}{L^\i}$, and $\norm{(w_r,w_z)}{L^\i}$ together in a bootstrap argument. Under the normalization $\norm{\bl{v}_0}{H^s}\le1$, we combine an $L^\i$ estimate for $v_r$ with the transport equation for $w_\th$ to obtain
$$
\int_0^t\norm{w_\th(t^\prime,\cd)}{L^\i}\d t^\prime\le C(1+t)^{7/3}
$$

throughout the bootstrap interval. We then combine this integral bound with a direct estimate of
$$
\int_0^t\norm{\na\bl{v}(t^\prime,\cd)}{L^\i}\d t^\prime
$$
to obtain a lifespan lower bound of order $[\log\log(\ve^{-1})]^{3/7}$ after normalization. The condition $s>7/2$, together with the Sobolev embedding
\[
H^s(\bR^3)\hookrightarrow W^{2,\infty}(\bR^3)\,,
\] 
provides the required $L^\i$ bounds for $\O_0$, $J_0$ and $\na(v_{0,r}/r)$. No compact support, sign condition, or spatial moment assumption is imposed.

For $\bl{v}_0\in H^3(\R^3)$, we retain the exponential bound for $\norm{w_\th(t,\cd)}{L^\i}$ obtained by Li and Zhou \cite{LZ}. By estimating $\int_0^t\|\bm{v}(t^\prime,\cd)\|_{L^\i}\d t^\prime$ directly, we remove one logarithm from their lifespan lower bound, improving the four-logarithm estimate to a three-logarithm estimate under the same $H^3$ regularity assumption.

\subsection{Statement of the main result}

\begin{theorem}\l{thm:main}
Fix $s>7/2$, and let $\bl{v}_0\in H^s(\R^3)$ be axisymmetric and divergence free, and $T_{*}$ be the maximal lifespan of the corresponding
$H^s$ solution of \ef{eq:euler}. Set
$$
\norm{\na\times (v_{0,\th}\etheta)}{L^\i}=\ve\ll 1.
$$
Then there exists a constant $C_s>0$, such that
\be\l{eq:H4-bound}
T_{*}\ge\f{C_s}{1+\norm{\bl{v}_0}{H^s}}
\lt[\log\log(\ve^{-1})\rt]^{3/7}.
\ee
\end{theorem}\qed

\begin{remark}\l{rem:three-log}
Let $\bl{v}_0\in H^3(\R^3)$ be axisymmetric and divergence free (same as in \cite{LZ}). Then there exits a constant $C>0$ such that
\be
\l{eq:new-bound}
T_{*}\ge\f{C}{1+\norm{\bl{v}_0}{H^3}}
\log\log\log({\ve^{-1}}).
\ee
Compared with \cite{LZ}, this estimate removes one iterated logarithm without requiring higher regularity of the initial data.

The proof is given at the end of the paper. It retains the weighted vorticity estimates of \cite{LZ} and closes the argument by estimating the accumulated velocity gradient $\int_0^t\norm{\na\bl{v}(t^\prime,\cd)}{L^\i}\d t^{\prime}$ directly.
\end{remark}\qed

\begin{remark}
We explain where the logarithms in \ef{eq:H4-bound}
and \ef{eq:new-bound} come from.
All estimates below are for normalized initial data
and under the a priori assumptions used in the proofs.

\begin{itemize}
\item
The first logarithm in \ef{eq:H4-bound}, farthest from
$\ve^{-1}$, comes from the estimate of
$\int_0^t\norm{\na\bl{v}(t^\prime,\cd)}{L^\i}\d t^{\prime}$.
Since $\norm{\na\bl{v}}{L^\i}$ cannot be controlled
by $\norm{\bl{w}}{L^\i}$ alone, we need a logarithmic
inequality. Combining it with the energy estimate
and the bound
$$
\int_0^t\norm{w_\th(t^\prime,\cd)}{L^\i}\d t^{\prime}\le C(1+t)^{7/3},
$$
Gronwall's inequality gives
$$
\int_0^t\norm{\na\bl{v}(t^\prime,\cd)}{L^\i}\d t^{\prime}
\le C\exp\lt(C(1+t)^{7/3}\rt).
$$

\item
When closing the continuity argument, we estimate
$(w_r,w_z)$ by applying the $L^\i$ transport estimate
and Gronwall's inequality:
\bn
\norm{(w_r,w_z)(t,\cd)}{L^\i}
&\le\ve\exp\lt(
C\int_0^t\norm{\na\bl{v}(t^\prime,\cd)}{L^\i}\d t^{\prime}
\rt)\\
&\le\ve\exp\lt[
C\exp\lt(C(1+t)^{7/3}\rt)
\rt].
\en
This gives the second logarithm in \ef{eq:H4-bound},
closest to $\ve^{-1}$.
The estimates of $\O$ and $J$
produce no further logarithm in the lifespan estimate.
\end{itemize}

For the $H^3$ result, we retain the estimate of $w_\th$ in \cite[Section~3.1, (3.7)]{LZ}. The improvement concerns the middle two logarithms described in \cite[Remark~1.2]{LZ}. In that proof, the full $H^3$ bound is used to control $\norm{\na\bl{v}}{L^\i}$ when estimating $(w_r,w_z)$; see \cite[(3.11)--(3.12)]{LZ}. We instead estimate the accumulated velocity gradient directly and obtain
\bn
\int_0^t\norm{\na\bl{v}(t^\prime,\cd)}{L^\i}\d t^{\prime}
&\le C\exp\lt(Ce^{Ct}\rt),\\
\norm{(w_r,w_z)(t,\cd)}{L^\i}
&\le\ve\exp\lt[C\exp\lt(Ce^{Ct}\rt)\rt].
\en
Unlike the approach in \cite{LZ}, here we do not use the full $H^3$ norm of $\bm{v}$ to bound $(w_r,w_z)$ in the bootstrap argument. This removes the third logarithm (counted from the outside) in \cite[(1.4)]{LZ} while the innermost logarithm coming from the transport estimate remains.
\end{remark}\qed

\subsection{Notations}
All spatial norms are taken over $\R^3$. For $1\le p\le\i$, we write $\norm{\cd}{L^p}$ for the
$L^p(\R^3)$ norm and use $\norm{\cd}{W^{s,p}}$ and $\norm{\cd}{\dot W^{s,p}}$ for the usual $L^p$-based inhomogeneous Sobolev norm and homogeneous Sobolev seminorm, respectively. When $p=2$, we use the standard notation
\[
H^s=W^{s,2},\qquad \dot H^s=\dot W^{s,2}.
\] Vector norms are defined using the Euclidean norm.

The symbol $\na^j$ denotes all Cartesian derivatives of order $j$, and $\p^\al$ denotes the derivative associated with a multi-index $\al$. For two operators $P,Q$, we write $[P,Q]=PQ-QP$. We also define
\be\l{eq:accumulated-gradient}
A(t)\ed\int_0^t\norm{\na\bl{v}(t^\prime,\cd)}{L^\i}\d t^{\prime}
\ee
the accumulated velocity gradient. The letter $C$ denotes a positive constant that may change from line to line. The notation $C_{m,M,\ldots}$ denotes a positive constant depending only on $m,M,\ldots$, which may also change from line to line.

\section{Preliminaries}

We first derive the equations used below. Taking the curl of the momentum equation in \ef{eq:euler},
we obtain
$$
\p_t\bl{w}+(\bl{v}\cd\na)\bl{w}
=(\bl{w}\cd\na)\bl{v}.
$$
Writing in the cylindrical coordinates $\bl{w}=w_r\er+w_\th\etheta+w_z\ez$, it follows that
\be
\l{eq:vorticity-system}
\begin{cases}
(\p_t+\bl{b}\cd\na)w_r
=(w_r\p_r+w_z\p_z)v_r,\\[6pt]
(\p_t+\bl{b}\cd\na)w_\th
=\f{v_r}{r}w_\th-2\f{v_\th}{r}w_r,\\[6pt]
(\p_t+\bl{b}\cd\na)w_z
=(w_r\p_r+w_z\p_z)v_z.
\end{cases}
\ee
Using $(\p_t+\bl{b}\cd\na)(\f{1}{r})=-\f{v_r}{r^2}$, we deduce from \ef{eq:vorticity-system} and the equation $\ef{eq:cylindrical}_2$ that
\be
\l{eq:weighted-system}
\begin{cases}
(\p_t+\bl{b}\cd\na)\O
=-2\f{v_\th}{r}J,\\[6pt]
(\p_t+\bl{b}\cd\na)J
=(w_r\p_r+w_z\p_z)\f{v_r}{r},\\[6pt]
(\p_t+\bl{b}\cd\na)\f{v_\th}{r}
=-2\f{v_r}{r}\f{v_\th}{r}.
\end{cases}
\ee

To estimate the vorticity components and the weighted quantities in \ef{eq:vorticity-system} and \ef{eq:weighted-system}, we recall the following standard $L^p$ estimate for transport equations with a divergence-free velocity field, Here we omit the proof, more details readers can refer to \cite{LZ}.
\begin{lemma}\l{lem:transport}
Let $1\le p\le\i$, $m\in\mathbb{N}$, and $\bl{b}\in L^1_{\mathrm{loc}}([0,T);W^{1,\i}(\R^3))$ be divergence free. Suppose $\bl{f}\in L^\i_{\mathrm{loc}}([0,T);L^p(\R^3;\R^m))$, with initial data $\bl{f}_0\in L^p$, is a weak solution of
$$
\p_t\bl{f}+\bl{b}\cd\na\bl{f}
=\bl{Q}\bl{f}+\bl{h}.
$$
Here $\bl{Q}(t,\cd):\R^3\to\R^{m\times m}$ and
$\bl{h}(t,\cd):\R^3\to\R^m$ satisfy
$$
\bl{Q}\in L^1_{\mathrm{loc}}([0,T);L^\i),\qq
\bl{h}\in L^1_{\mathrm{loc}}([0,T);L^p).
$$
Matrix norms are induced by the Euclidean norm. Then, for $0\le t<T$,
\bn
\norm{\bl{f}(t,\cd)}{L^p}
&\le\norm{\bl{f}_0}{L^p}
+\int_0^t\lt(
\norm{\bl{Q}(t^\prime,\cd)}{L^\i}\norm{\bl{f}(t^\prime,\cd)}{L^p}
+\norm{\bl{h}(t^\prime,\cd)}{L^p}\rt)\d t^{\prime},
\en
and
\bn
\norm{\bl{f}(t,\cd)}{L^p}
&\le\lt(\norm{\bl{f}_0}{L^p}
+\int_0^t\norm{\bl{h}(t^\prime,\cd)}{L^p}\d t^{\prime}\rt)
\exp\lt(\int_0^t\norm{\bl{Q}(t^\prime,\cd)}{L^\i}\d t^{\prime}\rt).
\en

\end{lemma}\qed

\begin{corollary}
Let $\bl{v}\in C([0,T);H^3(\R^3))$ be an axisymmetric solution of \ef{eq:euler}. Then, for $0\le t<T$,
\be
\l{eq:swirl-vorticity-transport}
\norm{(w_r,w_z)(t,\cd)}{L^\i}
\le\norm{(w_{0,r},w_{0,z})}{L^\i}e^{CA(t)},
\ee
where $A(t)$ is defined in \ef{eq:accumulated-gradient}.
\end{corollary}

\begin{proof}
Applying Lemma~\ref{lem:transport} with $p=\i$ to the first
and third equations of \ef{eq:vorticity-system}, whose
right-hand side coefficients are bounded by $C|\na\bl{v}|$,
we obtain \ef{eq:swirl-vorticity-transport}.
\end{proof}

The following axisymmetric Biot--Savart estimate is stated as Lemma~2.2 in \cite{LZ}; see also \cite[Proposition~2.5]{MZ}. We also omit the proof here.
\begin{lemma}\l{lem:axisymmetric}
Let $\bl{v}\in H^3(\R^3)$ be axisymmetric and divergence free. For $1<p<\i$, there exists an absolute constant $C_p>0$ such that
\be
\l{eq:axisymmetric-bs}
\norm{\na(\f{v_r}{r})(t,\cd)}{L^p}\le C_p\norm{\O(t,\cd)}{L^p}.
\ee
In particular, whenever $\O\in L^2\cap L^6$,
\be
\l{eq:radial-quotient-infty}
\norm{\f{v_r}{r}(t,\cd)}{L^\i}
\le C\norm{\O(t,\cd)}{L^2}^{1/2}\norm{\O(t,\cd)}{L^6}^{1/2}.
\ee
\end{lemma}\qed

We next establish estimates for the weighted quantities in \ef{eq:weighted-system}.
\begin{lemma}\l{lem:quotients}
Let $\bl{v}\in H^3(\R^3)$ be axisymmetric. Then
\be
\l{eq:swirl-quotient}
\norm{\f{v_\th}{r}(t,\cd)}{L^\i}\le\f12\norm{w_z(t,\cd)}{L^\i}.
\ee
Moreover,
\be
\l{eq:initial-weighted}
\norm{\O(t,\cd)}{L^2}+\norm{J(t,\cd)}{L^2}+\norm{\O(t,\cd)}{L^6}+\norm{J(t,\cd)}{L^6}
\le C\norm{\bl{v}}{H^3}.
\ee

If $\bl{v}\in H^s$ for some $s>7/2$, then
\be
\l{eq:axis-H4}
\norm{\O(t,\cd)}{L^\i}+\norm{J(t,\cd)}{L^\i}
+\norm{\na(\f{v_r}{r})(t,\cd)}{L^\i}\le C{\color{blue}_s}\norm{\bl{v}(t,\cd)}{H^s}.
\ee
\end{lemma}
\begin{proof}
The identity $\p_r(rv_\th)=rw_z$ and regularity at $r=0$ yield
$$
\f{v_\th(r,z)}r
=\f1{r^2}\int_0^r\r\,w_z(\r,z)\d\r,
$$
which proves \ef{eq:swirl-quotient}. 

To prove \ef{eq:initial-weighted}, recall that
$w_r,w_\th,w_z$ are independent of $\th$. Using $\p_\th\er=\etheta$, $\p_\th\etheta=-\er$,
and $\p_\th\ez=0$, we obtain
$$
\f1r\p_\th \bl{w}=J\etheta-\O\er.
$$
Since $r^{-1}\p_\th=\etheta\cd\na$ and $|\etheta|=1$, the preceding identity implies
$$
\lt|\O\rt|+\lt|J\rt|
\le \s{2}\,|\etheta\cd\na\bl{w}|
\le C|\na\bl{w}|,\qquad r>0.
$$
Using $\bl{w}=\na\times\bl{v}$ and the Sobolev inequality, we obtain
\bn
\norm{\O(t,\cd)}{L^2}
+\norm{J(t,\cd)}{L^2}
&\le C\norm{\na\bl{w}(t,\cd)}{L^2}
\le C\|\na^2\bl{v}(t,\cd)\|_{L^2},\\
\norm{\O(t,\cd)}{L^6}
+\norm{J(t,\cd)}{L^6}
&\le C\norm{\na\bl{w}(t,\cd)}{L^6}
\le C\|\na^2\bl{v}(t,\cd)\|_{L^6}
\le C\|\na^3\bl{v}(t,\cd)\|_{L^2}.
\en
Both right-hand sides are bounded by $C\norm{\bl{v}(t,\cd)}{H^3}$, which proves \ef{eq:initial-weighted}.

If $\bl{v}\in H^{(7/2)_+}$, the same pointwise estimate and the Sobolev embedding $H^{(7/2)_+}(\R^3)\hookrightarrow W^{2,\i}(\R^3)$ give
\bn
\norm{\O(t,\cd)}{L^\i}
+\norm{J(t,\cd)}{L^\i}
&\le C\norm{\na\bl{w}(t,\cd)}{L^\i}\le C\|\na^2\bl{v}(t,\cd)\|_{L^\i}
\le C_s \norm{\bl{v}(t,\cd)}{H^s},\q\text{for some}\q s>\f{7}{2}\,.
\en
This proves the bounds for the first two terms in \ef{eq:axis-H4}. The last term is similar.
\end{proof}

We now recall conservation of kinetic energy and a standard Sobolev estimate for solutions of \ef{eq:euler}.
\begin{lemma}\l{lem:energy}
Let $s>5/2$, and let $\bl{v}$ be an $H^s$ solution of
\ef{eq:euler}. Then
\be
\l{eq:energy}
\norm{\bl{v}(t,\cd)}{L^2}
=\norm{\bl{v}_0}{L^2},\qq
\norm{\bl{v}(t,\cd)}{H^s}
\le \norm{\bl{v}_0}{H^s}
\exp\lt(
C_s\int_0^t
\norm{\na\bl{v}(\tau,\cd)}{L^\i}\,\mathrm d\tau
\rt).
\ee
\end{lemma}

\begin{proof}
We first carry out the calculation for smooth solutions.
Taking the dot product of the momentum equation with
$\bl{v}$, integrating over $\R^3$, and using
$\div\bl{v}=0$, we obtain
\bn
\f12\f{\mathrm d}{\mathrm dt}
\norm{\bl{v}(t,\cd)}{L^2}^2
&=-\int_{\R^3}
(\bl{v}\cd\na)\bl{v}\cd\bl{v}\dx
-\int_{\R^3}\na p\cd\bl{v}\dx\\
&=\f12\int_{\R^3}
(\div\bl{v})|\bl{v}|^2\dx
+\int_{\R^3}p\,\div\bl{v}\dx\\
&=0.
\en
Integrating in time gives
\[
\norm{\bl{v}(t,\cd)}{L^2}
=\norm{\bl{v}_0}{L^2}.
\]

We next estimate the fractional derivatives of order $s$.
Let
\[
\Lambda=(1-\Delta)^{1/2},
\]
so that
\[
\norm{\Lambda^s\bl{v}}{L^2}
=\norm{\bl{v}}{H^s}.
\]
Applying $\Lambda^s$ to the momentum equation gives
\[
\p_t\Lambda^s\bl{v}
+(\bl{v}\cd\na)\Lambda^s\bl{v}
+\na\Lambda^s p
=-[\Lambda^s,\bl{v}\cd\na]\bl{v}.
\]
Taking the dot product with $\Lambda^s\bl{v}$ and
integrating over $\R^3$, we use
\[
\div\bl{v}=0,
\qquad
\div(\Lambda^s\bl{v})
=\Lambda^s\div\bl{v}=0,
\]
to eliminate the transport and pressure terms. Indeed,
\[
\int_{\R^3}
(\bl{v}\cd\na)\Lambda^s\bl{v}
\cd\Lambda^s\bl{v}\dx
=-\f12\int_{\R^3}
(\div\bl{v})|\Lambda^s\bl{v}|^2\dx=0,
\]
and
\[
\int_{\R^3}\na\Lambda^s p
\cd\Lambda^s\bl{v}\dx
=-\int_{\R^3}\Lambda^s p\,
\div(\Lambda^s\bl{v})\dx=0.
\]
Consequently,
\be
\l{eq:fractional-energy}
\f12\f{\mathrm d}{\mathrm dt}
\norm{\Lambda^s\bl{v}}{L^2}^2
=-\int_{\R^3}
[\Lambda^s,\bl{v}\cd\na]\bl{v}
\cd\Lambda^s\bl{v}\dx.
\ee

Writing
\[
[\Lambda^s,\bl{v}\cd\na]\bl{v}
=\sum_{j=1}^3[\Lambda^s,v_j]\p_j\bl{v},
\]
the Kato--Ponce commutator estimate \cite{KP} gives
\bn
\|[\Lambda^s,\bl{v}\cd\na]\bl{v}\|_{L^2}
&\le C_s\sum_{j=1}^3
\lt(
\norm{\na v_j}{L^\i}
\|\Lambda^{s-1}\p_j\bl{v}\|_{L^2}
+\norm{\Lambda^s v_j}{L^2}
\norm{\p_j\bl{v}}{L^\i}
\rt)\\
&\le C_s\norm{\na\bl{v}}{L^\i}
\norm{\bl{v}}{H^s}.
\en
It follows from \ef{eq:fractional-energy} and the
Cauchy--Schwarz inequality that
\bn
\f12\f{\mathrm d}{\mathrm dt}
\norm{\bl{v}}{H^s}^2
&\le
\norm{[\Lambda^s,\bl{v}\cd\na]\bl{v}}{L^2}
\norm{\Lambda^s\bl{v}}{L^2}\\
&\le
C_s\norm{\na\bl{v}}{L^\i}
\norm{\bl{v}}{H^s}^2.
\en
Gronwall's inequality therefore yields
\[
\norm{\bl{v}(t,\cd)}{H^s}^2
\le
\norm{\bl{v}_0}{H^s}^2
\exp\lt(
C_s\int_0^t
\norm{\na\bl{v}(\tau,\cd)}{L^\i}\,\mathrm d\tau
\rt).
\]
After taking the square root and changing the value of
$C_s$, we obtain
\[
\norm{\bl{v}(t,\cd)}{H^s}
\le
\norm{\bl{v}_0}{H^s}
\exp\lt(
C_s\int_0^t
\norm{\na\bl{v}(\tau,\cd)}{L^\i}\,\mathrm d\tau
\rt).
\]

Since $s>5/2$, Sobolev embedding gives
$H^s(\R^3)\hookrightarrow W^{1,\i}(\R^3)$, so all the quantities
appearing in the estimate are controlled for an $H^s$
solution. By smooth approximation and local well-posedness
in $H^s$, the same estimates hold for the stated solutions
on every compact time interval within their lifespan.
This proves \ef{eq:energy}.
\end{proof}

The following logarithmic inequality is a consequence of
\cite[Lemma~2.5]{LZ}; see also the references therein.
We combine it with the preceding energy estimate to control the accumulated velocity gradient $A(t)$.
\begin{lemma}\l{lem:log}
For a divergence-free vector field $\bl{g}\in H^3(\R^3)$,
\be
\l{eq:log-gradient}
\norm{\na \bl{g}}{L^\i}
\le C\lt[1+\norm{\bl{g}}{L^2}+
\norm{\na\times \bl{g}}{L^\i}\log\bigl(e+\norm{\bl{g}}{H^3}\bigr)\rt].
\ee
In particular, if $\norm{\bl{g}}{L^2}\le1$, the first two terms on the
right-hand side can be replaced by an absolute constant.

Moreover, if $\bl{v}$ solves \ef{eq:euler} and $\norm{\bl{v}_0}{H^s}\le1$,
$s\geq3$, then
\be
\l{eq:A-differential}
A'(t)\le C\lt[1+
\bigl(\norm{w_\th(t,\cd)}{L^\i}+\norm{(w_r,w_z)(t,\cd)}{L^\i}\bigr)
(1+A(t))\rt].
\ee
\end{lemma}
\begin{proof}
We only prove \ef{eq:A-differential} here. Under the normalization $\norm{\bl{v}_0}{H^s}\le1$,
\ef{eq:energy} gives
$$
\norm{\bl{v}(t,\cd)}{L^2}\le1,\qq
\norm{\bl{v}(t,\cd)}{H^3}
\le\norm{\bl{v}(t,\cd)}{H^s}\le Ce^{CA(t)}.
$$
Hence
$$
\log\bigl(e+\norm{\bl{v}(t,\cd)}{H^3}\bigr)
\le C(1+A(t)).
$$
By the definition of $A(t)$ and \ef{eq:log-gradient}, we obtain
\bn
A'(t,\cd)=\norm{\na\bl{v}(t,\cd)}{L^\i}
&\le C\lt[1+\norm{\bl{w}(t,\cd)}{L^\i}
\log\bigl(e+\norm{\bl{v}(t,\cd)}{H^3}\bigr)\rt]\\
&\le C\lt[1+
\bigl(\norm{w_\th(t,\cd)}{L^\i}
+\norm{(w_r,w_z)(t,\cd)}{L^\i}\bigr)(1+A(t))\rt],
\en
which proves \ef{eq:A-differential}.
\end{proof}

The next estimate follows from local kernel bounds in \cite[Section~3.1]{LJ} and a far-field $L^2$ estimate in \cite[Appendix~A.1]{LJ}. No compact support or moment assumptions are needed. The proof uses the axisymmetry of $\bl{b}$ and the identities $\div\bl{b}=0$ and $\na\times\bl{b}=w_\th\etheta$, without using an evolution equation for $\bl{b}$.

\begin{lemma}[Radial velocity bound]\l{lem:radial-global}
Let $\bl{b}=v_r\er+v_z\ez\in H^s(\R^3)$, $s>7/2$ be axisymmetric and
divergence free, with $\na\times\bl{b}=w_\th\etheta$ and
$\O=w_\th/r$. Suppose that
$$
\norm{\bl{b}}{L^2}\le1,\qq
\norm{\O}{L^\i}\le K,\qq K\ge1.
$$
Then
\be
\l{eq:radial-global}
\norm{v_r}{L^\i}
\le C_K\lt(1+\norm{w_\th}{L^\i}^{1/4}\rt),
\ee
where $C_K$ depends only on $K$.
\end{lemma}

\begin{proof}
We first obtain a local estimate at $(1,0)$. Use $(\r,\zeta)\in(0,\i)\times\R$ as the integration variables and write
$$
D=((1-\r)^2+\zeta^2)^{1/2}
$$
for the distance to $(1,0)$ in the meridional half-plane. With our convention for vorticity, the axisymmetric Biot-Savart formula (see \cite[Section~2.2]{LJ}) reads
\bn
v_r(1,0)
&=\int_0^\i\int_\R
\mK(\r,\zeta)(\na\times\bl{b})_\th(\r,\zeta)\d\zeta\d\r,\\
\mK(\r,\zeta)
&=-\f{\r\zeta}{4\pi}\int_{-\pi}^{\pi}
\f{\cos\phi\d\phi}
{(1+\r^2+\zeta^2-2\r\cos\phi)^{3/2}}.
\en
Here, the axisymmetric Biot–Savart kernel given in \cite[(2.2)]{LJ} is evaluated at the observation point $(1,0)$. Choose a smooth function $\chi:[0,\i)\to[0,1]$ with $\chi=1$ on $[0,1]$ and $\chi=0$ on $[2,\i)$, and split the integral using $1=\chi(4D)+(1-\chi(4D))$.

The far-field $L^2$ estimate in \cite[Appendix~A.1]{LJ} gives
\be
\l{eq:radial-far-energy}
\lt|\int_0^\i\int_\R
(1-\chi(4D))\mK(\na\times\bl{b})_\th\d\zeta\d\r\rt|
\le C\norm{\bl{b}}{L^2}.
\ee
Indeed, integration by parts and the Cauchy--Schwarz inequality give this estimate, as in the bounds for $(I''_{11}),\ldots,(I''_{14})$ in \cite[Appendix~A.1]{LJ}, before the interpolation leading to their (A.5). The kernel bounds in their (A.2)--(A.4) ensure that the boundary term at the axis vanishes and that the cutoff error at infinity tends to zero. We omit the details. 

For the near-field part, \cite[Lemma~2.1 and Section~3.1]{LJ} gives
\be
\l{eq:kernel-near}
|\mK|\le CD^{-1},\qq
|\na_{\r,\zeta}\mK|\le CD^{-2},\qq 0<D<1/2.
\ee
Let $0<\ell\le1/8$. Noticing that $\chi(4D)\chi(D/\ell)=\chi(D/\ell)$, we use the decomposition
$$
\chi(4D)=\chi(D/\ell)
+\chi(4D)(1-\chi(D/\ell)).
$$
First we see
\bn
\lt|\int_0^\i\int_\R
\chi(D/\ell)\mK(\na\times\bl{b})_\th\d\zeta\d\r\rt|
&\le C\sup_{D\le2\ell}
\lt|(\na\times\bl{b})_\th(\r,\zeta)\rt|
\int_{D\le2\ell}D^{-1}\d\r\d\zeta\\
&\le C\ell\sup_{D\le2\ell}
\lt|(\na\times\bl{b})_\th(\r,\zeta)\rt|.
\en
On the other hand, the truncated kernel $\chi(4D)(1-\chi(D/\ell))\mK$ is supported in $\ell\le D\le1/2$, where $1/2\le\r\le3/2$. The derivatives of $\chi(4D)$ and $\chi(D/\ell)$ are supported in $1/4\le D\le1/2$ and $\ell\le D\le2\ell$, respectively. Consequently,
$$
|\na_{\r,\zeta}\chi(4D)|
+|\na_{\r,\zeta}\chi(D/\ell)|
\le CD^{-1}.
$$
Together with \ef{eq:kernel-near}, the product rule gives
$$
\lt|
\na_{\r,\zeta}
\lt[\chi(4D)(1-\chi(D/\ell))\mK\rt]
\rt|
\le CD^{-2}
\qquad\text{for }\ell\le D\le1/2.
$$
Since $\r$ is bounded above and below on this region, integration by parts and the Cauchy--Schwarz inequality give
\bn
&\lt|\int_0^\i\int_\R
\chi(4D)(1-\chi(D/\ell))\mK
(\na\times\bl{b})_\th\d\zeta\d\r\rt|\\
&\q\le C\norm{\bl{b}}{L^2}
\lt(\int_{\ell\le D\le1/2}
D^{-4}\d\r\d\zeta\rt)^{1/2}
\le C\ell^{-1}\norm{\bl{b}}{L^2}.
\en

Combining these estimates with \ef{eq:radial-far-energy}, we obtain
\be\l{Evrr}
|v_r(1,0)|\le C\lt[
(1+\ell^{-1})\norm{\bl{b}}{L^2}
+\ell\sup_{D\le2\ell}
\lt|(\na\times\bl{b})_\th(\r,\zeta)\rt|\rt].
\ee
This estimate holds for every $0<\ell\le1/8$ and does not use
the bounds $\norm{\bl{b}}{L^2}\le1$ or $\norm{\O}{L^\i}\le K$.

We now fix $r\ge1$ and $z\in\R$, and set
$$
\wt{\bl{b}}(x)=\bl{b}(rx+z\ez).
$$
Then
\bn
\wt v_r(1,0)&=v_r(r,z),\q\|\wt{\bl{b}}\|_{L^2}
&=r^{-3/2}\norm{\bl{b}}{L^2}\le r^{-3/2},\q (\na\times\wt{\bl{b}})_\th(\r,\zeta)
&=r w_\th(r\r,z+r\zeta).
\en
On the region $D\le2\ell\le1/4$, we have
$$
|\r-1|\le D\le\f14,\qq \f12\le\r\le\f32.
$$
Since $w_\th=r\O$, it follows that
$$
|w_\th(r\r,z+r\zeta)|
\le Kr\r\le2Kr.
$$
Combining this with
$|w_\th(r\r,z+r\zeta)|\le\norm{w_\th}{L^\i}$ gives
\bn
\sup_{D\le2\ell}
\lt|(\na\times\wt{\bl{b}})_\th(\r,\zeta)\rt|
&\le r\min\lt\{\norm{w_\th}{L^\i},2Kr\rt\}\le r^{3/2}\lt(2K\norm{w_\th}{L^\i}\rt)^{1/2}.
\en
Applying \ef{Evrr} to $\wt{\bl{b}}$, we obtain
$$
|v_r(r,z)|\le C\lt[
(1+\ell^{-1})r^{-3/2}
+\ell r^{3/2}\lt(K\norm{w_\th}{L^\i}\rt)^{1/2}\rt].
$$
Choose
$$
\ell=\lt[8+r^{3/2}\lt(K\norm{w_\th}{L^\i}\rt)^{1/4}\rt]^{-1}.
$$
Then $0<\ell\le1/8$, and
\bn
(1+\ell^{-1})r^{-3/2}
&=9r^{-3/2}+\lt(K\norm{w_\th}{L^\i}\rt)^{1/4},\\
\ell r^{3/2}\lt(K\norm{w_\th}{L^\i}\rt)^{1/2}
&=\f{r^{3/2}\lt(K\norm{w_\th}{L^\i}\rt)^{1/2}}
{8+r^{3/2}\lt(K\norm{w_\th}{L^\i}\rt)^{1/4}}
\le\lt(K\norm{w_\th}{L^\i}\rt)^{1/4}.
\en
Consequently,
\be\l{Efar1018}
|v_r(r,z)|\le C\lt(1+K^{1/4}\norm{w_\th}{L^\i}^{1/4}\rt),
\qq r\ge1.
\ee

It remains to consider $r\le1$. Fix $x$ with radial
coordinate at most $1$. The three-dimensional Biot-Savart formula is
$$
\bl{b}(x)=\f1{4\pi}\int_{\R^3}
\f{(\na\times\bl{b})(y)\times(x-y)}
{|x-y|^3}\dy.
$$
Using the same cutoff $\chi$, split the integral with
$\chi(|x-y|)$.
In the far-field term, integrate by parts to transfer the curl
from $\bl{b}$ to the kernel. The differentiated kernel and the
terms containing derivatives of $\chi$ are bounded by
$C|x-y|^{-3}$ on $|x-y|\ge1$.
This integration by parts is justified in $L^2$, since the
far-field kernel and its first derivatives belong to $L^2$.
Consequently,
\bn
|\bl{b}(x)|
&\le C\int_{|x-y|\le2}
\f{|(\na\times\bl{b})(y)|}{|x-y|^2}\dy
+C\int_{|x-y|\ge1}
\f{|\bl{b}(y)|}{|x-y|^3}\dy\\
&\le C\norm{\na\times\bl{b}}{L^\i(B(x,2))}
+C\norm{\bl{b}}{L^2}
\lt(\int_{|x-y|\ge1}
|x-y|^{-6}\dy\rt)^{1/2}.
\en
Here $B(x,2)$ denotes the ball of radius $2$ centered at
$x$. Every point in this ball has radial coordinate at most $3$,
so $|\na\times\bl{b}|=|w_\th|\le3K$ there.
It follows that
\be\l{ENear1018}
|v_r(r,z)|\le|\bl{b}(x)|\le C(K+1),\qq r\le1.
\ee
Combining estimates \eqref{Efar1018} and \eqref{ENear1018}, we prove \ef{eq:radial-global}.
\end{proof}

\section{Proof of the main theorem}

\begin{proof}[Proof of Theorem~\ref{thm:main}]
We first normalize the initial norm. Set $M=1+\norm{\bl{v}_0}{H^s}$ and define
\be
\l{eq:scaling}
\wt{\bl{v}}(t,x)=M^{-1}\bl{v}(t/M,x),\qq
\wt{p}(t,x)=M^{-2}p(t/M,x).
\ee
The rescaled fields solve \ef{eq:euler}, the initial $H^s$ norm is at most $1$, and the initial swirl vorticity has norm $\ve/M$. We omit the tildes and denote the normalized swirl size again by $\ve$. Thus, throughout this proof,
$$
\norm{\bl{v}_0}{H^s}\le1,\qq
\norm{(w_{0,r},w_{0,z})}{L^\i}=\ve.
$$
By \ef{eq:axis-H4}, we can fix an absolute constant $K\ge1$ such that
$$
\norm{\O_0}{L^\i}\le K,\qq
\norm{J_0}{L^\i}\le K.
$$
We impose the bootstrap assumptions
\be
\l{eq:H4-bootstrap}
\norm{(w_r,w_z)(t,\cd)}{L^\i}\le1,\qq
\norm{\O(t,\cd)}{L^\i}\le2K,\qq
\norm{J(t,\cd)}{L^\i}\le2K.
\ee
They hold on a nonempty time interval for some $\ve\ll 1$.

\medskip\noindent\emph{Time integral of the azimuthal vorticity.}
The second equation of \ef{eq:vorticity-system} reads
\be\l{Ewth1018}
(\p_t+\bl{b}\cd\na)w_\th
=v_r\O-2\f{v_\th}{r}w_r.
\ee
Energy conservation gives
$\norm{\bl{b}(t,\cd)}{L^2}\le\norm{\bl{v}(t,\cd)}{L^2}\le1$.
Since $\na\times\bl{b}=w_\th\etheta$ and $\bl{b}\in H^s$,
Lemma~\ref{lem:radial-global}, with $2K$ in place of $K$, and
\ef{eq:H4-bootstrap} imply
$$
\norm{v_r(t,\cd)}{L^\i}
\le C_K\lt(1+\norm{w_\th(t,\cd)}{L^\i}^{1/4}\rt).
$$
Moreover, \ef{eq:swirl-quotient} and \ef{eq:H4-bootstrap} give
\be\l{eq:swirl-source-bound}
2\norm{\f{v_\th}{r}(t,\cd)}{L^\i}
\norm{w_r(t,\cd)}{L^\i}
\le\norm{(w_r,w_z)(t,\cd)}{L^\i}^2\le1.
\ee
Combining the radial velocity bound with \ef{eq:swirl-source-bound}
and \ef{eq:H4-bootstrap}, we obtain
\ba\l{eq:theta-source-global}
\norm{\lt(v_r\O-2\f{v_\th}{r}w_r\rt)(t,\cd)}{L^\i}
&\le\norm{v_r(t,\cd)}{L^\i}\norm{\O(t,\cd)}{L^\i}+2\norm{\f{v_\th}{r}(t,\cd)}{L^\i}\norm{w_r(t,\cd)}{L^\i}\\
&\le C_K\lt(1+\norm{w_\th(t,\cd)}{L^\i}^{1/4}\rt).
\ea

Thus, applying Lemma~\ref{lem:transport} with $p=\i$, $m=3$, and
$\bl{Q}=0$, and using \ef{eq:theta-source-global}, we obtain from \eqref{Ewth1018} that
\ba\l{wthinf1018}
\norm{w_\th(t,\cd)}{L^\i}
&\le\norm{w_{0,\th}}{L^\i}
+\int_0^t\norm{\lt(v_r\O-2\f{v_\th}{r}w_r\rt)(t^\prime,\cd)}{L^\i}\d t^{\prime}\\
&\le C_K(1+t)+C_K\int_0^t
\norm{w_\th(t^\prime,\cd)}{L^\i}^{1/4}\d t^{\prime}.
\ea
Here we used the normalized initial bound and Sobolev embedding.
Set
$$
F(t)\ed\int_0^t\norm{w_\th(t^\prime,\cd)}{L^\i}\d t^{\prime}.
$$
By H\"older's inequality, one can deduce
$$
\int_0^t\norm{w_\th(t^\prime,\cd)}{L^\i}^{1/4}\d t^{\prime}
\le t^{3/4}
\lt(\int_0^t\norm{w_\th(t^\prime,\cd)}{L^\i}\d t^{\prime}\rt)^{1/4}
=t^{3/4}F(t)^{1/4}.
$$
Then one infers from \eqref{wthinf1018} that
\bn
F'(t)&=\norm{w_\th(t,\cd)}{L^\i}\le C_K(1+t)+C_Kt^{3/4}F(t)^{1/4}.
\en
Since $F(0)=0$ and $F$ is nondecreasing,
integration in time gives
\bn
F(t)
&\le C_K(t+t^2)+C_K\int_0^t {(t^{\prime})}^{3/4} F(t^\prime)^{1/4}\d t^{\prime}\\
&\le C_K(t+t^2)+C_KF(t)^{1/4}\int_0^t ({t^\prime})^{3/4}\d t^{\prime}\\
&\le C_K(t+t^2)+C_Kt^{7/4}F(t)^{1/4}.
\en
By Young's inequality with exponents $4$ and $4/3$, one deduces
$$
C_Kt^{7/4}F(t)^{1/4}
\le\f12F(t)+C_Kt^{7/3}.
$$
Absorbing $\f12F(t)$ into the left-hand side and using
$t+t^2\le (1+t)^{7/3}$, we obtain
\be
\l{eq:theta-integral}
\int_0^t\norm{w_\th(t^\prime,\cd)}{L^\i}\d t^{\prime}
=F(t)\le C(1+t)^{7/3}.
\ee
The initial vorticity norm and $K$ are bounded by absolute constants
under the normalization $\norm{\bl{v}_0}{H^s}\le1$.
Thus $C$ is independent of the initial data and $\ve$.

\medskip\noindent\emph{High-order estimates.}
Lemma~\ref{lem:energy} gives
$$
\norm{\bl{v}(t,\cd)}{H^s}\le Ce^{CA(t)},\qq
\log\bigl(e+\norm{\bl{v}(t,\cd)}{H^s}\bigr)\le C(1+A(t)).
$$
The common estimate \ef{eq:A-differential} therefore applies. It follows from
\ef{eq:H4-bootstrap} and \ef{eq:theta-integral} that
\bn
1+A(t)
&\le\exp\lt[Ct+C\int_0^t
\bigl(\norm{w_\th(t^\prime,\cd)}{L^\i}+\norm{(w_r,w_z)(t^\prime,\cd)}{L^\i}\bigr)
\d t^{\prime}\rt]\\
&\le\exp\lt(Ct+CF(t)\rt)\\
&\le C\exp\lt(C(1+t)^{7/3}\rt).
\en
Combining \ef{eq:energy} and \ef{eq:swirl-vorticity-transport} with this bound gives
\ba
\l{eq:H4-growth}
\norm{\bl{v}(t,\cd)}{H^s}
&\le C\exp\lt[C\exp\lt(C(1+t)^{7/3}\rt)\rt],
\ea
\bn
\norm{(w_r,w_z)(t,\cd)}{L^\i}
&\le\ve\exp\lt[C\exp\lt(C(1+t)^{7/3}\rt)\rt].
\en
In particular, multiplication does not add an exponential level:
$$
\norm{(w_r,w_z)(t,\cd)}{L^\i}\norm{\bl{v}(t,\cd)}{H^s}
\le C\ve
\exp\lt[C\exp\lt(C(1+t)^{7/3}\rt)\rt].
$$

\medskip\noindent\emph{Closing the weighted-vorticity bounds.}
The equations for $J$ and $\O$ in \ef{eq:weighted-system}, together
with \ef{eq:axis-H4} and \ef{eq:swirl-quotient}, imply
\ba
\l{eq:radial-weighted-H4}
\norm{J(t,\cd)}{L^\i}
&\le\norm{J_0}{L^\i}+C\int_0^t
\norm{(w_r,w_z)(t^\prime,\cd)}{L^\i}\norm{\bl{v}(t^\prime,\cd)}{H^s}\d t^{\prime},
\ea
\ba
\l{eq:theta-weighted-H4}
\norm{\O(t,\cd)}{L^\i}
&\le\norm{\O_0}{L^\i}+
\int_0^t\norm{(w_r,w_z)(t^\prime,\cd)}{L^\i}
\norm{J(t^\prime,\cd)}{L^\i}\d t^{\prime},\\
&\le\norm{\O_0}{L^\i}+
2K\int_0^t\norm{(w_r,w_z)(t^\prime,\cd)}{L^\i}\d t^{\prime}.
\ea
Thus $\norm{(w_r,w_z)(t,\cd)}{L^\i}$ and both increments in
\ef{eq:radial-weighted-H4}--\ef{eq:theta-weighted-H4} are bounded by
\be
\l{eq:H4-error}
C\ve(1+t)
\exp\lt[C\exp\lt(C(1+t)^{7/3}\rt)\rt].
\ee

Set
$$
T_\ve=c[\log\log(\ve^{-1})]^{3/7}.
$$
First choose $c$ sufficiently small, and then choose $\ve$
sufficiently small, so that
$$
C(1+T_\ve)^{7/3}\le\f12\log\log(\ve^{-1}).
$$
The logarithm of \ef{eq:H4-error} at $T_\ve$ is at most
$$
-\log(\ve^{-1})+
\log\bigl(C(1+T_\ve)\bigr)+
C\bigl[\log(\ve^{-1})\bigr]^{1/2},
$$
which tends to $-\i$ as $\ve\searrow 0$, this indicate \ef{eq:H4-error} tends to zero. For a fixed absolute smallness threshold, we therefore obtain
$$
\norm{(w_r,w_z)(t,\cd)}{L^\i}\le\f12,\qq
\norm{\O(t,\cd)}{L^\i}\le\f32K,\qq
\norm{J(t,\cd)}{L^\i}\le\f32K,
\qq 0\le t\le T_\ve,
$$
throughout the bootstrap interval. These inequalities strictly improve \ef{eq:H4-bootstrap}. The three bootstrap quantities are continuous by the $H^s$ local theory and \ef{eq:axis-H4}. Thus none can first reach its threshold before $\min\{T_\ve,T_{*}\}$. If $T_{*}\le T_\ve$, \ef{eq:H4-growth} gives a uniform $H^s$ bound up to $T_{*}$, contradicting the continuation theorem. The normalized lifespan is therefore at least $T_\ve$. Restoring the original time scale and the original value of $\ve$ proves \ef{eq:H4-bound}.
\end{proof}

\section{Proof of the three-logarithm improvement}

\begin{proof}[Proof of Remark~\ref{rem:three-log}]
Apply \ef{eq:scaling} with $M=1+\norm{\bl{v}_0}{H^3}$ and use the same convention for the normalized variables and for $\ve$. We therefore assume
\be
\l{eq:normalized}
\norm{\bl{v}_0}{H^3}\le1,\qq
\norm{(w_{0,r},w_{0,z})}{L^\i}=\ve.
\ee
Estimate \ef{eq:swirl-vorticity-transport} reads
$$
\norm{(w_r,w_z)(t,\cd)}{L^\i}\le\ve e^{CA(t)}.
$$
Local existence and continuity provide a nonempty interval on which
\be
\l{eq:bootstrap}
t\sup_{0\le t^{\prime}\le t}\norm{(w_r,w_z)(t^\prime,\cd)}{L^\i}\le1.
\ee
We first establish the estimates below on any interval satisfying
\ef{eq:bootstrap}. In particular,
\be
\l{eq:integrated-swirl-bound}
\int_0^t\norm{(w_r,w_z)(t^\prime,\cd)}{L^\i}\d t^{\prime}\le1.
\ee

\medskip\noindent\emph{Estimate from Li and Zhou.}

The condition \ef{eq:bootstrap} is the a priori assumption (3.2) of Li and Zhou \cite{LZ}, since
$$
\norm{\na\times(v_\th\etheta)}{L^\i}
=\norm{(w_r,w_z)}{L^\i}.
$$
We use their estimate \cite[Section~3.1, (3.7)]{LZ}, which follows from the weighted-vorticity equations (3.3) and the estimates (3.4)--(3.6) in that paper. Under the normalization
\ef{eq:normalized}, it gives
\be
\l{eq:theta-exponential}
\norm{w_\th(t,\cd)}{L^\i}\le Ce^{Ct}.
\ee
Thus we retain the vorticity estimate of Li and Zhou and turn to a direct bound for the accumulated velocity gradient.

\medskip\noindent\emph{Higher-order bounds and the continuity argument.}

Lemma~\ref{lem:energy} and \ef{eq:normalized} give
\be
\l{eq:sobolev-through-gradient}
\norm{\bl{v}(t,\cd)}{H^3}\le Ce^{CA(t)},\qq
\log\bigl(e+\norm{\bl{v}(t,\cd)}{H^3}\bigr)\le C(1+A(t)).
\ee
Applying Gronwall's inequality to \ef{eq:A-differential}, and then using
\ef{eq:integrated-swirl-bound} and \ef{eq:theta-exponential}, we obtain
\ba
\l{eq:A-bound}
1+A(t,\cd)
&\le\exp\lt[Ct+C\int_0^t
\bigl(\norm{w_\th(t^\prime,\cd)}{L^\i}+\norm{(w_r,w_z)(t^\prime,\cd)}{L^\i}\bigr)
\d t^{\prime}\rt]\\
&\le C\exp(Ce^{Ct}).
\ea
Thus the $H^3$ norm stays finite on every finite interval satisfying
\ef{eq:bootstrap}.

We now close the continuity argument using \ef{eq:swirl-vorticity-transport} directly. There exist absolute constants $C_1,C_2,C_3>0$ such that
\be
\l{eq:swirl-vorticity-bound}
\sup_{0\le t^{\prime}\le t}\norm{(w_r,w_z)(t^\prime,\cd)}{L^\i}
\le\ve\exp\lt[C_1\exp\lt(C_2e^{C_3t}\rt)\rt].
\ee
Choose $c>0$ so that $C_3c<1$, and let
$$
T_\ve=c\log\log\log(\ve^{-1}).
$$
For $0\le t\le T_\ve$,
$$
C_1\exp(C_2e^{C_3t})
\le C_1\exp\lt[C_2\bigl(\log\log(\ve^{-1})\bigr)^{C_3c}
\rt]
=o\bigl(\log(\ve^{-1})\bigr)
\q\text{as }\ve\searrow0.
$$
After choosing an absolute smallness threshold for $\ve$, we therefore have
$$
C_1\exp(C_2e^{C_3T_\ve})
\le\f12\log(\ve^{-1}),\qq
T_\ve\ve^{1/2}\le\f12.
$$
It follows from \ef{eq:swirl-vorticity-bound} that
\be
\l{eq:bootstrap-improved}
\norm{(w_r,w_z)(t,\cd)}{L^\i}\le\ve^{1/2},\qq
t\sup_{0\le t^\prime\le t}\norm{(w_r,w_z)(t^\prime,\cd)}{L^\i}\le\f12,
\qq 0\le t\le T_\ve,
\ee
as long as the solution exists and the bootstrap condition holds.

For completeness, let $T_b$ be the supremum of the times below the maximal existence time for which the strict inequality in \ef{eq:bootstrap} holds. This number is positive by local existence and continuity. If $T_b<T_\ve$ and $T_b$ is below the maximal existence time, \ef{eq:bootstrap-improved} and continuity extend the bootstrap interval past $T_b$, a contradiction. If the maximal existence time were at most $T_\ve$,
\ef{eq:A-bound} and \ef{eq:sobolev-through-gradient} would give a uniform
$H^3$ bound up to that finite time. The continuation theorem for Euler would extend the solution, again a contradiction. Thus the normalized lifespan is at least $T_\ve$.

Restoring the original time scale and the original value of $\ve$ in \ef{eq:scaling}, we conclude that
$$
T_{*}\ge\f{C}{1+\norm{\bl{v}_0}{H^3}}
\log\log\log\lt({\ve}^{-1}\rt).
$$
The a priori calculations can first be performed for smooth axisymmetric data and then passed to $H^3$ data by axisymmetry-preserving approximation and local well-posedness. The norms used in the bootstrap are continuous in time. This proves \ef{eq:new-bound}.
\end{proof}

\end{document}